\documentclass[12pt]{article}
\usepackage{amsmath,amsthm,amssymb}
\usepackage{mathtools}

\newtheorem{theorem}{Theorem}
\newtheorem{lemma}{Lemma}
\newtheorem{cor}{Corollary}

\numberwithin{equation}{section}

\newcommand{\D}{{\rm d}}
\newcommand{\dx}{\, \D x}

\newcommand{\ds}{\, \D s}
\newcommand{\dr}{\, \D r}
\newcommand{\dt}{\, \D t}

\newcommand{\rz}{\mathbb{R}}

\newcommand{\eps}{\varepsilon}
\newcommand{\iom}{\int_{\Omega}}

\newcommand{\klauf}{\left(\begin{array}}
\newcommand{\klzu}{\end{array}\right)}

\title{An extension of a theorem of Bers and Finn on the removability of isolated singularities to the
Euler-Lagrange equations related to general linear growth problems}

\author{Michael Bildhauer \& Martin Fuchs}
\date{}

\newcommand{\reff}[1]{(\ref{#1})}

\usepackage{hyperref}

\hypersetup{
   pdfnewwindow=true,
   colorlinks=false,
   linkcolor=black,
   citecolor=green,
   filecolor=magenta,
   urlcolor=cyan,
}

\begin{document}

\parindent0em
\maketitle

\newcommand{\op}[1]{\operatorname{#1}}
\newcommand{\bv}{\op{BV}}
\newcommand{\mub}{\overline{\mu}}
\newcommand{\muhat}{\hat{\mu}}

\newcommand{\hypref}[2]{\hyperref[#2]{#1 \ref*{#2}}}
\newcommand{\hypreff}[1]{\hyperref[#1]{(\ref*{#1})}}

\newcommand{\ob}[1]{^{(#1)}}

\newcommand{\xh}{\Xi}
\newcommand{\oh}[1]{O\left(#1\right)}
\newcommand{\xn}{\hat{x}}
\newcommand{\yn}{\hat{y}}
\newcommand{\On}{\hat{\Omega}}

\newcommand{\hn}{\hat{N}}

\begin{abstract}{}
A famous theorem of Bers and Finn states that isolated singularities of solutions to the non-parametric
minimal surface equation are removable. We show that this result remains valid, if the area functional is replaced
by a general functional of linear growth depending on the modulus of the gradient.\\

We emphasize that Serrin (\cite{Se:1965_1}) in fact proved the removability of singularities on sets of $(n-1)$-dimensional Hausdorff
measure zero in an even more general setting. Our main interest is to generalize the comparison principles as outlined, for instance,
in Section 10 of  \cite{Os:1986_1} without having the particular geometric structure of minimal surfaces. 
It turns out that generalized catenoids serve as an appropriate tool for proving our results.
\end{abstract}

\parindent0ex
\section{Introduction}\label{intro}

We discuss solutions $u \in C^2(D)$ defined on an open set $D \subset \rz^n$ of the equation
\begin{equation}\label{intro 1}
\op{div} \Bigg[\frac{g'\big(|\nabla u|\big)}{|\nabla u|} \nabla u \Bigg] = 0
\end{equation}

arising as the Euler-Lagrange equation of the variational problem
\begin{equation}\label{intro 2}
\iom g \big(|\nabla v|\big) \dx \to \min
\end{equation}

among functions $v$: $\Omega \to \rz$ with prescribed boundary data. The assumptions concerning the density $g$
are as follows:\\

we consider functions $g$: $[0,\infty) \to \rz$ of class $C^{2,\alpha}\big([0,\infty)\big)$ for some exponent $0 < \alpha < 1$
being of linear growth in the sense
that with suitable constants $a$, $A >0$, $b$, $B \geq 0$ the inequality
\begin{equation}\label{intro 3}
at - b \leq g(t) \leq A t + B
\end{equation}

holds for any $t \geq 0$. Moreover, we require strict convexity of $g$ by imposing the condition
\begin{equation}\label{intro 4}
g''(t) > 0  \quad\mbox{for all}\quad  t \geq 0 \, .
\end{equation}

Finally, we assume
\begin{equation}\label{intro 5}
g'(0) = 0 \, .
\end{equation}

We then will prove that the famous theorem of Bers and Finn (see \cite{Be:1951_1}, \cite{Fi:1953_1})
on the removability of isolated singularities for solutions of the non-parametric
minimal surface equation extends to any solution of \reff{intro 1} provided that $g$ satisfies these hypotheses.\\

In more detail we have the following result:

\begin{theorem}\label{intro theo 1}
Consider an open set $\Omega \subset \rz^n$, fix some point $x_0 \in \Omega$ and assume that
$u \in C^2\big(\Omega -\{x_0\}\big)$ is a solution of equation \reff{intro 1}
on the set $D:= \Omega -\{x_0\}$ with $g$ satisfying \reff{intro 3}--\reff{intro 5}.

Then $u$ admits an extension $\overline{u}\in C^2(\Omega)$ and $\overline{u}$ solves equation \reff{intro 1} on the set $\Omega$.
\end{theorem}

In the case of minimal surfaces, i.e.~for the choice $g(t) = \sqrt{1+t^2}$ in equation \reff{intro 1} and for $n=2$,
the result of the theorem was proved 
independently by Bers \cite{Be:1951_1} and Finn \cite{Fi:1953_1}. 
Concerning solutions of the non-parametric minimal surface equation in dimensions $n>2$
the removability of singular sets $K$ being closed subsets of $\Omega$ 
such that $\mathcal{H}^{n-1}(K) = 0$ was established by 
DeGiorgi and Stampacchia \cite{DeGS:1965_1}, by Simon \cite{Si:1977_4}, 
Anzellotti \cite{An:1981_1} and Miranda \cite{Mi:1977_1}.\\

During the proof  of  Theorem \ref{intro theo 1} we will have to distinguish two essentially different cases,
where the first one is closely related to the minimal surface setting in the sense that we suppose
\begin{equation}\label{intro 6}
\int_0^\infty t g''(t) \dt < \infty
\end{equation}
restricting the growth of $g''$ at infinity. Note that \reff{intro 6} is a consequence of the pointwise inequality
\begin{equation}\label{intro 7}
g''(t) \leq c (1+t)^{-\mu}\, , \quad t \geq 0 \, ,
\end{equation}

provided we choose $\mu > 2$. In the minimal surface case, i.e.~for the choice $g(t) = \sqrt{1+t^2}$, we
can choose $\mu =3$ in estimate \reff{intro 7}, and by
a ``$\mu$-surface in $\rz^{n+1}$'' we denote the graph $\big\{\big(x,u(x)\big) \in \rz^{n+1}:\, x \in D\big\}$ of a solution
$u$: $D \to \rz$ of equation \reff{intro 1}, provided that $g$ satisfies the conditions \reff{intro 3}, \reff{intro 4}
and \reff{intro 7} for some exponent $\mu > 2$. We refer to the recent manuscript \cite{BF:2021_2} on some geometric
properties of $\mu$-surfaces in the case $n=2$. Adopting this notation we deduce from Theorem \ref{intro theo 1} that $\mu$-surfaces
do not admit isolated singular points.\\ 

However, this removability property does not depend on any geometric features.  
As it is formulated in Theorem \ref{intro theo 1}, the non-existence of isolated singularities is just
a consequence of the linear growth of $g$ which is also exploited
in the second case
\begin{equation}\label{intro 8}
\int_0^\infty t g''(t) \dt = \infty \, .
\end{equation}
This condition already occurs, e.g., in \cite{BBM:2018_1} (compare also \cite{BF:2020_2}) in a quite different setting: 
in Theorem 1.1 of \cite{BBM:2018_1}, equation
\reff{intro 8} together with some kind of balancing condition serves as a criterion for the solvability of a classical Dirichlet-problem,
where the authors argue with the help of suitable barrier functions. Both in \cite{BBM:2018_1} and in \cite{BF:2020_2}
generalized catenoids are used as basic tools, which is also the case in our considerations. 
Depending on the conditions \reff{intro 6} and \reff{intro 8}, respectively,
these catenoids are of infinte height or uniformly bounded.

\section{Proof of Theorem \ref{intro theo 1} under condition \reff{intro 6}}\label{proof 1}

In the following we consider energy densities $g$: $[0,\infty) \to \rz$ of class $C^{2,\alpha}$ such that \reff{intro 3}--\reff{intro 5} hold.
In particular $g'$ is a bounded function and strictly increasing, thus
\begin{equation}\label{s1 1}
0 = g'(0) < g'(t) \to :  g'_\infty \quad\mbox{as}\quad t \to \infty \, .
\end{equation}

W.l.o.g.~it is assumed that
\begin{equation}\label{s2 2}
g'_\infty = 1 \, .
\end{equation}

Moreover, $g \in C^2$ together with $g'(0) = 0$ yields that the function $G$: $\rz^n \to \rz$,
$G(p) := g\big(|p|\big)$ is of class $C^2(\rz^n)$ satisfying
\begin{equation}\label{s1 3}
\sum_{i,j =1}^n \frac{\partial^2 G}{\partial p_i \partial p_j}(p) q_i q_j > 0 \quad\mbox{for all}\quad p,\, q\in \rz^n \, , \quad q \not=0 \, .
\end{equation}

\vspace*{2ex}
\emph{Step 1. Maximum principle.}\\

We observe that in the subsequent considerations we may not assume \mbox{Lipschitz} continuity of solutions up to the boundary, hence
Theorem 1.2 of \cite{Mi:1997_1} does not apply. We will make use of the following variant:

\begin{lemma}\label{s2 lem 1}
Suppose that $D$ is a bounded Lipschitz domain in $\rz^n$ and that we have \reff{intro 3}--\reff{intro 5}. Moreover suppose
that $u$, $v \in C^2(D) \cap C^0(\overline{D})$ satisfy equation \reff{intro 1}. Then we have:
\[
u \leq v+M \quad\mbox{on $\partial D$ for some real number $M$} \quad \Rightarrow\quad
u \leq v+ M\quad\mbox{on $D$}\, .
\]
\end{lemma}

\vspace*{2ex}
With Lemma \ref{s2 lem 1} the following corollary is immediate:
\begin{cor}\label{s2 cor 1}
The Dirichlet-problem associated to \reff{intro 1} within the class $C^2(D) \cap C^0(\overline{D})$ admits at most one solution.
\end{cor}

\vspace*{2ex}
\emph{Proof of Lemma \ref{s2 lem 1}.} From \reff{intro 1} one obtains
\begin {equation}\label{s2 1}
0 = \sum_{i,j=1}^n \frac{\partial^2 G}{\partial p_i \partial p_j}(\nabla u) \partial_{x_i} \partial_{x_j} u \quad\mbox{on $D$}\ .
\end{equation}

Now we refer to Theorem 10.1, p.~263, of \cite{GT:1998_1} with coefficients 
($x \in D$, $y\in \rz$, $p \in \rz^n$)
\[
a_{ij}(x,y,p) := \frac{\partial^2 G}{\partial p_i \partial p_j} (p)
\]

which, by \reff{s1 3}, are seen to be elliptic. Since we consider the admissible function space $C^2(D) \cap C^0(\overline{D})$ the proof
is complete with the above mentioned reference. \hfill\qed\\

\emph{Step 2. Generalized catenoids as comparison surfaces.}\\

Let $g$ satisfy \reff{intro 3}--\reff{intro 5} and \reff{intro 6} and recall \reff{s1 1}, which implies that $g'$ maps $[0,\infty)$ in a one-to-one way
onto the intervall $[0,1)$.\\ 

For numbers $\alpha >0$ and constants $a \in \rz$ we define for $x\in \rz^n$, $ |x| >\alpha^{1/(n-1)}$,
\begin{eqnarray}\label{s3 1}
k_{\alpha, a}^{\pm}(x)&:=& l_{\alpha ,a}^{\pm}\big(|x|\big)\nonumber\\
&:=& \pm \int_{\alpha^{1/(n-1)}}^{|x|} \big(g'\big)^{-1}\Big(\frac{\alpha}{r^{n-1}}\Big) \dr + a  \, .
\end{eqnarray}

We have:

\begin{lemma}\label{s3 lem 1}
The functions $k_{\alpha , a}^{\pm}$ are solutions of problem \reff{intro 1} on $|x| > \alpha^{1/(n-1)}$ with continuous extension (through
the value $a$)  to the boundary $|x| = \alpha^{1/(n-1)}$.
\end{lemma}

\vspace*{2ex}
\emph{Proof of Lemma \ref{s3 lem 1}.} W.l.o.g. we let $\alpha =1$, $a=0$ in the definition of $k_{\alpha ,a}^{\pm}$ and $l_{\alpha ,a}^{\pm}$, respectively.
Dropping the indices $\alpha$ and $a$ we have for $t>1$ (letting $r^{n-1}=1/g'(s)$)
\begin{eqnarray*}
l^+(t) &=& \int_1^t \big(g'\big)^{-1}\Big(\frac{1}{r^{n-1}}\Big) \dr\\
&=& \int_\infty^{s^*(t)} s g''(s) \Bigg[-\frac{1}{n-1} \big(g'(s)\big)^{-\frac{n}{n-1}}\Bigg] \ds \, ,\\
s^*&=&s^*(t):= \big(g'\big)^{-1}\Big(\frac{1}{t^{n-1}}\Big) \, ,
\end{eqnarray*}

hence
\begin{equation}\label{s3 2}
l^+(t) = \int_{s^*(t)}^\infty s g''(s)\Bigg[\frac{1}{n-1} \big(g'(s)\big)^{-\frac{n}{n-1}}\Bigg]  \ds 
\end{equation}
is well defined at least for $t>1$ on account of assumption \reff{intro 6} and due to the behaviour of $g'$ as stated  in \reff{s1 1}.\\

Moreover, from \reff{s3 2} it immediately follows that 
\[
\lim_{t\downarrow 1}l^+(t) = 0 \, ,
\]
thus $l^+$ has a continuous extension to $t=1$ by letting $l^+(1) = 0$.\\

Let us look at equation \reff{intro 1} in the case
\[
D=B_R(0) - \overline{B_r(0)}
\]
for balls centered at $0$ with radii $0 < r < R \leq \infty$. Suppose further that we have a solution $u(x)$ of the form $u(x) = \varphi(\rho)$,
$\rho = |x|$. Then \reff{intro 1} is equivalent to the ODE
\begin{equation}\label{s3 3}
\frac{\D}{\D \rho} \Bigg[\rho^{n-1} \frac{g'\big(|\varphi'(\rho)|\big)}{|\varphi'(\rho)|} \varphi'(\rho) \Bigg] = 0\, , \quad
\rho\in (r,R)\, ,
\end{equation}
and obviously $l^+$ solves \reff{s3 3} for the choices $r=1$, $R=\infty$. This proves Lemma \ref{s3 lem 1},
since with obvious modifications the above calculations can be adjusted to the functions
$k^{-}_{\alpha,a}$. \qed\\

\emph{Step 3. Comparison principle.}\\

\begin{lemma}\label{s4 lem 1}Let $g$ satisfy the assumptions \reff{intro 3} -- \reff{intro 6}. For $0 < r < R < \infty$ let
$D = B_R(0) -\overline{B_r(0)}$ in equation \reff{intro 1} and consider a solution $u \in C^2(D) \cap C^1\big(\overline{D}\big)$
such that for some $a\in \rz$ it holds with $\alpha := r^{n-1}$
\begin{equation}\label{s4 1}
u \leq k^{-}_{\alpha ,a} \quad\mbox{on}\quad \partial B_R(0) \, .
\end{equation}
Then we have
\begin{equation}\label{s4 2}
u \leq k^{-}_{\alpha,a} \quad\mbox{throughout}\quad\overline{D}\, .
\end{equation}
\end{lemma}

\emph{Proof of Lemma \ref{s4 lem 1}.} By Lemma \ref{s2 lem 1} and Lemma \ref{s3 lem 1} it is enough to show that
\begin{equation}\label{s4 3}
u \leq k^{-}_{\alpha ,a} \quad{on}\quad \partial B_r(0) \, .
\end{equation}

W.l.o.g. let $r=1$ and $a=0$ and write $k^-$ in place of $k^{-}_{1,a}$.
Following a standard reasoning known from the minimal surface case (compare
\cite{Os:1986_1}) we assume that \reff{s4 3} is wrong. Then we can choose $x_0 \in \partial B_1 (0)$ satisfying (on account of 
$k^- \equiv 0$ on $\partial B_1(0)$)
\begin{equation}\label{s4 4}
0 < u(x_0) = \max_{|x|=1} u(x) =: M\, .
\end{equation}
For $t > 1$ we let
\[
\phi(t) := u\big(t x_0\big) - k^{-}\big(t x_0\big)
\]
and get
\[
\phi'(t) = x_0 \cdot \nabla u\big(t x_0\big) + \big(g'\big)^{-1} \Big(\frac{1}{t^{n-1}}\Big) \, .
\]
Since we assume $u \in C^1\big(\overline{D}\big)$ and since we have 
\[
\big(g'\big)^{-1} \Big(\frac{1}{t}\Big) \to \infty \quad\mbox{as}\quad t \downarrow 1\, ,
\]
there exists $\eps >0$ such that $\phi'(t) > 0$ for all $t \in (1,1+\eps)$.
This implies
\begin{equation}\label{s4 5}
u\big(tx_0\big) - k^{-}\big(tx_0\big) > u(x_0) \quad\mbox{on}\quad (1,1+\eps) \, .
\end{equation}
Recalling the definition of $M$ and our assumption \reff{s4 1}, Lemma \ref{s2 lem 1} yields
\begin{equation}\label{s4 6}
u -k^{-} \leq M \quad\mbox{on}\quad \overline{D} \, .
\end{equation}
Obviously \reff{s4 6} contradicts \reff{s4 5}, thus we have \reff{s4 3} and the proof is complete. \qed\\

\emph{Step 4. Removability of isolated singularities.}\\

Now we are going to prove the first part of the theorem. Let $g$ satisfy \reff{intro 3} -- \reff{intro 6} and consider
a solution $u \in C^2\big(B_R(0)-\{0\}\big)$ of equation \reff{intro 1} on the punctured ball
$B_R(0) - \{0\}$. W.l.o.g. we assume that $u \in C^1\big(\overline{B_R(0)} - B_r(0)\big)$ for any radius $0 < r < R$. 
Following standard arguments (compare \cite{Os:1986_1}) we claim 
\begin{equation}\label{s5 1}
\min_{|x|=R} u \leq u(y) \leq \max_{|x| =R} u(x) \quad \mbox{for all}\quad y \not= 0\, ,\quad |y| < R\, .
\end{equation}
In fact, we let 
\[
M(r) := \max_{|x|=r} u(x) \, , \quad 0 < r < R \, ,
\]
and define 
\[
a:= M(R) + \int_r^R \big(g'\big)^{-1}\Big(\frac{\alpha}{t^{n-1}}\Big) \dt \, ,
\]
i.e.~we have (again with $\alpha = r^{n-1})$
\begin{equation}\label{s5 2}
k^{-}_{\alpha ,a} \equiv M(R)\quad\mbox{on}\quad \partial B_R(0) \, , \quad\mbox{hence}\quad
u \leq k^{-}_{\alpha ,a}\quad\mbox{on}\quad \partial B_R(0) \, .
\end{equation}
Quoting Lemma \ref{s4 lem 1} and observing that \reff{s5 2} corresponds to hypothesis \reff{s4 1}, we obtain (compare \reff{s4 2})
\begin{equation}\label{s5 3}
u \leq k^{-}_{\alpha,a}\quad \mbox{on}\quad \overline{B_R(0)}-B_r(0)\, .
\end{equation}
Fix a point $x$ such that $ 0 < |x| < R$. Then \reff{s5 3} implies for any $0 < r = \alpha^{1/(n-1)} < |x|$ (recall \reff{s3 1} and \reff{s5 2})
\begin{equation}\label{s5 4}
u(x) \leq k^{-}_{\alpha ,a}(x) = M(R) + \int_{|x|}^{R} \big(g'\big)^{-1}\Big(\frac{\alpha}{t^{n-1}}\Big) \dt \, .
\end{equation}
Recall that $x$ is fixed and that we have \reff{intro 5}. Hence, passing to the limit $r \to 0$ in \reff{s5 4}, we obtain
\[
u(x) \leq M(R)  \quad\mbox{for any}\quad x \in B_R(0) -\{0\} \, ,
\]
thus the second inequality stated in \reff{s5 1} is established. The first inequality in \reff{s5 1} follows with obvious modifications.\\

Finally, let $\tilde{u} \in C^1\big(\overline{B_R(0)}\big)$ be a smooth extension of  $u_{|\partial B_{R}(0)}$.
From the Hilbert-Haar theory we find a unique Lipschitz-minimizer $v$: $\overline{B_R(0)} \to \rz$ of the energy
\[
\int_{B_R(0)} G\big(\nabla w\big) \dx = \int_{B_R(0)} g \big(|\nabla w|\big) \dx
\]
subject to the boundary data $\tilde{u}_{|\partial B_{R}(0)} = u_{|\partial B_{R}(0)}$, which due to our hypotheses 
(recall that $g \in C^{2,\alpha}\big([0,\infty)\big)$) turns out to be 
of class $C^{2,\beta}\big(B_R(0)\big)$ for some $\beta \in (0,1)$.
In fact, the Hilbert-Haar minimizer $v$ has H\"older continuous first derivatives (see, e.g.~\cite{Mi:1997_1}, Theorem 1.7)
and standard arguments from regularity theory applied to equation \reff{s2 1} imply $v\in C^{2,\beta}$.\\

We claim $v=u$ on $B_R(0)-\{0\}$, which means that $v$ is the desired $C^2$-extension of $u$.\\

In fact, for $0 < \eps \ll 1$ it holds (using \reff{s2 1} for the functions $u$ and $v$ on $D=B_R(0)-\{0\}$
and with $\nu$ denoting the exterior normal on
$\partial (B_R(0) - B_\eps(0))$)
\begin{eqnarray*}
\lefteqn{ \int_{B_R(0) - B_\eps(0)} (\nabla u- \nabla v) \big(DG(\nabla u) - DG(\nabla v)\big) \dx}\\
&=& \int_{\partial (B_R(0) - B_\eps(0))} (u-v) \big(DG(\nabla u) - DG(\nabla v)\big) \cdot \nu\,  \D \mathcal{H}^{1}\\
& \to & 0 \quad \mbox{as}\quad \eps \to 0
\end{eqnarray*}
on account of $u$, $v \in L^{\infty}\big(B_R(0)\big)$. By ellipticity this implies $\nabla u = \nabla v$ on $B_R(0) - \{0\}$
and our claim follows. \qed\\

\section{Proof of Theorem \ref{intro theo 1} under condition \reff{intro 8}}\label{infin}
Let the density $g$ satisfy the same assumptions as stated in the beginning of the previous section, in particular we
have \reff{s2 2}, but now we replace \reff{intro 6} by condition \reff{intro 8}.
W.l.o.g.~we may assume that $\Omega = B_2(0)$, $x_0 = 0$ in Theorem \ref{intro theo 1}.
Replacing \reff{s3 1} we now fix $0 < r <1$ and let 
\begin{eqnarray}\label{fin 1}
k^{\pm}_{r,a}(x)&:=& l_{r ,a}^{\pm}\big(|x|\big)\nonumber\\
&:=&a  \pm \int_{1}^{|x|} \big(g'\big)^{-1}\Big(\frac{r}{t^{n-1}}\Big) \dt \, , \quad |x| > r^{\frac{1}{n-1}}\, .
\end{eqnarray}
Then we have for $|x|>r^{1/(n-1)}$ (letting $t^{n-1}=r/g'(s)$ for the fixed number $0 < r < 1$)
\begin{eqnarray*}
k_{r,a}^{\pm}(x) &=& a \mp r^{\frac{1}{n-1}} \int_{(g')^{-1}(r)}^{s^*(|x|)} s g''(s) \Bigg[\frac{1}{n-1} \big(g'(s)\big)^{-\frac{n}{n-1}}\Bigg] \ds \, ,\\
s^*&=&s^*(|x|):= \big(g'\big)^{-1}\Big(\frac{r}{|x|^{n-1}}\Big) \, .
\end{eqnarray*}
We recall \reff{intro 8} and note that $k^{\pm}_{r,a}(x)$ is defined for $|x| > r^{1/(n-1)}$ with limit
\begin{equation}\label{fin 2}
k^{\pm}_{r,a}(x) \to \mp \infty \quad\mbox{as}\quad |x| \to r^{\frac{1}{n-1}}\, .
\end{equation}
Here $a \in \rz$ is chosen according to (note $k^{\pm}_{r,a}(x) = a$ for $|x| =1$)
\begin{equation}\label{fin 3}
u \leq k^{-}_{r,a} = a := \max_{\partial B_1(0)} u \quad \mbox{on}\quad \partial B_1(0)\,  .
\end{equation}

Since $u$ is bounded on $\partial B_{r^{1/(n-1)}}(0)$ and since we have \reff{fin 2}, 
we may choose $\eps >0$ sufficiently small such that
\begin{equation}\label{fin 4}
u \leq k^{-}_{r,a}\quad\mbox{on}\quad \partial B_{(r+\eps)^{1/(n-1)}}(0) \, .
\end{equation}

With \reff{fin 3} and \reff{fin 4} we now directly apply Lemma \ref{s2 lem 1} to obtain
\begin{equation}\label{fin 5}
u \leq k^{-}_{r,a} \quad\mbox{on}\quad \overline{B_1(0)} - B_{(r+\eps)^{1/(n-1)}}(0) \, .
\end{equation}

W.l.o.g.~we may suppose $\eps < r/2$, hence \reff{fin 5} gives
\begin{equation}\label{fin 6}
u \leq k^{-}_{r,a} \quad\mbox{on}\quad \overline{B_1(0)} - B_{(3r/2)^{1/(n-1)}}(0) \, .
\end{equation}

Now we fix $x \in B_1(0)$, $x \not= 0$, and recall that the real number $a$ chosen in \reff{fin 3} is not depending on 
the radius $r$ considered above. We let
\begin{equation}\label{fin 7}
r:= \frac{1}{2} |x|^{n-1}\, , \quad\mbox{i.e.}\quad |x| = (2r)^{1/(n-1)} \, .
\end{equation}

Then \reff{fin 6} together with the choice \reff{fin 7} finally yields
\begin{equation}\label{fin 8}
u(x) \leq k^{-}_{r,a}\big(x) \, .
\end{equation}

On the other hand definition \reff{fin 1}
gives together with the monotonicity of $g'$
\begin{eqnarray}\label{fin 9}
k^{-}_{r,a}(x)=l^{-}_{r,a}\big(|x|\big) &=& a - \int_1^{|x|} \big(g')^{-1}\Big(\frac{1}{2} \frac{|x|^{n-1}}{t^{n-1}}\Big)\dt\nonumber\\
&=& a + \int_{|x|}^1 \big(g')^{-1}\Big(\frac{1}{2} \frac{|x|^{n-1}}{t^{n-1}}\Big)\dt\nonumber\\
&\leq & a + \big(1-|x|\big) \big(g'\big)^{-1} \big(1/2\big) \leq c \, ,
\end{eqnarray}
where the constant $c$ is not depending on $x$, hence we have established an uniform upper bound for $u$ and a uniform
lower bound follows along similar lines.
Proceeding exactly as done  in the first case at the end of Step 4 the theorem is proved. \qed


\bibliography{Bers_01_06_22}
\bibliographystyle{unsrt}

\end{document}